\documentclass[11pt]{article}
\usepackage{amsmath, amsthm,  amsfonts}
\input{amssym.def}

\def\R{{\mathbb R}}

\def\C{{\mathbb C}}
\def\H{{\mathbb H}}
\def\E{{\mathbb E}}
\def\S{{\mathbb S}}

\def\tilde{\widetilde}

\def\dis{\displaystyle}

\def\and{\hbox{\quad and \quad}}

\def\Bbb{\bf}
\catcode`\@=12
\newtheorem{theorem}{Theorem}
\newtheorem{definition}{Definition}

\theoremstyle{remark}

\newtheorem{remark}{Remark}

\input{epsf.tex}

\begin{document}

\vskip 0.3cm
\begin{center}
{\LARGE  {\bf The Hyperbolic Geometry of the Sinh-Gordon Equation}}

\

Preprint
\end{center}

\begin{center}
{\bf {Magdalena Toda\footnote{Department of Mathematics and Statistics, Texas Tech
University, Lubbock, Texas 79409-1042, U.S.A. {\em mtoda@math.ttu.edu}}\\}}
\end{center}

{\abstract This preliminary report studies immersed surfaces of constant mean
curvature in $H^3$ through their {\it adjusted Gauss maps} (as harmonic maps in $S^2$)
and their {\it adjusted frames} in $SU(2)$. Lawson's correspondence between Euclidean
CMC surfaces and their hyperbolic cousins is interpreted here under a different
perspective: the equivalence of their Weierstrass representations (normalized
potentials). This work also presents a construction algorithm for the moving frame,
the adjusted frame, their Maurer-Cartan forms, and ultimately the CMC immersion.}

\section{Introduction}

The classical Weierstrass representation formula for minimal surfaces in $\E^3$
consists of a meromorphic function (the classical Gauss map) and a holomorphic
$1$-form. Several years ago, a method now referred to as DPW (\cite{DoPeWu}) was
introduced for nonzero constant mean curvature surfaces (abbreviated as CMC) in $E^3$.
The method gives a characterization of these surfaces in terms of a pair of functions,
called normalized (or meromorphic) potential, and also a method to construct all
associate immersions based on loop group factorization. Among the classes of surfaces
for which such a potential was found explicitly, and used to construct surfaces, we
mention: constant mean curvature surfaces in Euclidean space $\E^3$ (\cite{DoHa});
minimal surfaces in $\E^3$ (\cite{DoPeTo}); weakly regular pseudospherical surfaces in
$\E^3$ (\cite{To}); timelike surfaces in the Minkowski 3-space (\cite{DoInTo});
Willmore surfaces in $\E^3$(\cite{He}); timelike minimal surfaces (\cite{TI}).

\

In this paper, we consider only cousins of the Euclidean CMC surfaces, that is,
surfaces with constant mean curvature $|H|> 1$. The case $|H|<1$, consisting of
solutions to the cosh-Gordon equation, has only one distinguished representative, the
minimal surfaces $H=0$. The author has devoted an entirely different approach (and
separate paper) to this case. More precisely, for minimal surfaces, the adjusted
Maurer-Cartan and the adjusted (harmonic) Gauss map are different, and they require
different loop groups.

As a byproduct, we here show that for any associated family of CMC surfaces with
$|H|>1$, there is a (strongly conformal) corresponding family, obtained by spectral
deformation, whose frame is unitary ($\rm SU(2)$-valued, as opposed to $\rm SL(2, \Bbb
C)$). This is not the case for $|H|<1$.

\

A first step of this study was communicated in 2002 as a preliminary report, but a DPW
method for constant mean surfaces in $\H^3$ was yet to be found. The author then
introduced the normalized potentials and conjectured that they represent the input for
a construction algorithm. A proposed DPW algorithm was publicized in 2004, through an
outstanding work by four authors, \cite{Ross}, for surfaces with constant mean
curvature $|H|>1$.

\

Our spectral deformations, potentials and DPW method are different from \cite{Ross}.
We study the normalized potential and show that it basically reduces to a `Weierstrass
pair': the Hopf differential, together with the holomorphic part of the metric
conformal factor. We here analyze how the normalized potential is used in order to
generate the {\it adjusted $\rm SU(2)$ frame}, the regular $\rm SL(2, \C)$ frame,
their two Maurer-Cartan forms, and the CMC immersion.

\ Among the advantages of this particular approach are the simple form of the
normalized potential compared to other representations one may use, as well as the
fact that one does not have to keep track of the monodromy representation. Also,
although the factorizations are not explicit, the resulting frames and immersions are.

\

We prompt the reader to check and note the following:

Even for the spectral deformation from \cite{Ross}, the usual frame $F$ is {\it not}
r-unitary. Due to the off-diagonal entries of the Lax matrices, containing $H-1$ and
$H+1$, the usual Maurer-Cartan form $F^{-1} \cdot dF$ is {\it not} ${\rm
su}(2)$-valued (even when the spectral parameter takes values on $S^1$).

The r-unitarization that is aimed in \cite{Ross} (see (1.10), (2.4) and (3.1))
actually takes place in our context, for a different reason (see formula (17)).

Let $\H^3(-1)$ denote the hyperbolic 3-space of constant sectional curvature $-1$.
Surfaces of constant mean curvature $|H|=1$ represented the topic of many papers over
the past fifteen years. An important result of R. Bryant (\cite{Br}) gave a
representation formula for these surfaces. In \cite{Ai},  R. Aiyama  and K. Akutagawa
gave a Kenmotsu-Bryant type representation formula for (branched) surfaces in
$\H^3(-c^2)$ of constant mean curvature $|H|\geq c$. In \cite{Ai2}, the same two
authors showed even further that there exists a Kenmotsu-Bryant type representation
formula for surfaces in $\H^3(-c^2)$ of constant mean curvature $|H|<c$.

\

A general result known as Lawson's correspondence has the following theorem as a
corollary:

\begin{theorem}
There is a bijective correspondence between the space of isometric immersions of
constant mean curvature $H>0$ in $\E^3$ and the space of isometric immersions of
constant mean curvature $\sqrt{H^2+1}$ in $\H^3(-1)$.
\end{theorem}

From now on, we will assume $\H^3$ as being of sectional curvature $-1$, unless
otherwise stated.

Through the Lawson correspondence, CMC surfaces in $\E^3$ are corresponded to CMC
$|H|>1$ surfaces in $\H^3$. On the other hand, minimal surfaces in $\E^3$ are
corresponded to CMC $|H|=1$ surfaces in $\H^3$. In the past decade, significant
progress has been made in the area of surfaces in $\H^3$, especially surfaces of
constant mean curvature greater than one. Recently, there has been some progress in
visualizing some surfaces. For example, N. Schmitt used a loop group splitting in
order to construct surfaces of constant mean curvature $|H|>1$ in $\H^3$. Based on
these methods, he wrote a program that produces hyperbolic analogues of some CMC
surfaces in $\E^3$, such as CMC bubbletons, CMC cylinders, Smyth surfaces, and
N-noids. The pictures of these CMC surfaces in $\H^3$ can be viewed at the GANG's
gallery of CMC surfaces (http://www.gang.umass.edu/).

\section{Integrable Systems of Constant Mean Curvature Surfaces in Hyperbolic
3-Space $\H^3$}

Let us consider the $4$-dimensional Lorentzian space
$$\R^{3,1}=\{(x^0,x^1,x^2,x^3)|
ds^2=-(dx^0)^2+(dx^1)^2+(dx^2)^2+(dx^3)^2\}.$$

The hyperbolic 3-space is the spacelike 3-manifold
$$\H^3=\H^3(-1)=\{x\in\R^{3,1}|<x,x>=-1, x^0>0\}$$
of constant sectional curvature $-1$.

Note that the following correspondence
$$x=(x^0,x^1,x^2,x^3)\longmapsto x=\left(\begin{array}{cc}
x^0+x^3 & x^1+ix^2\\
x^1-ix^2 & x^0-x^3
\end{array}\right)$$
provides an identification between $\R^{3,1}$ and the space of $2\times 2$ Hermitian
matrices. \ The complex Lie group ${\rm SL}(2,\C)$ acts isometrically and transitively
on $\H^3(-1)$ by
\begin{eqnarray*}
{\rm SL}(2,\C)\times\H^3(-1)&\buildrel\cdot\over\longrightarrow&\H^3(-1)\\
(g,h)&\longmapsto&g\cdot h=ghg^*,
\end{eqnarray*}
where $g^*={\bar{g}}^t$. Thus, $\H^3={\rm SL}(2,\C)/{\rm SU}(2)$.

\

Let $M$ be a simply connected Riemann surface and $f: M\longrightarrow\H^3$ an
immersion.

Consider $(e^0=f,e^1,e^2,e^3)$ the local orthonormal frame of the immersion $f$. Then
we have
\begin{eqnarray*}
&{}&de^0=df=\omega_ie^i,\ i=1,2,\\
&{}&de^j=\omega_je^0+\omega^j_ie^i,\ i=1,2,3,
\end{eqnarray*}
where $\omega^i_j=-\omega^j_i$ and $\omega^i_i=0$.

For the adapted frame of the immersion $f$, Cartan's structure equations can be
written on short as
\begin{eqnarray*}
d\omega_i=\omega^j_i\wedge\omega_j\\
d\omega^i_j+\omega^i_k\wedge\omega^k_j+\omega_i\wedge\omega_j=0.
\end{eqnarray*}
Let $\sigma_i,\ i=0,1,2,3$,  be the following matrices
$$\sigma_0=\left(\begin{array}{cc}
1 & 0\\
0 & 1
\end{array}\right),\ \sigma_1=\left(\begin{array}{cc}
0 & 1\\
1 & 0
\end{array}\right),\ \sigma_2=\left(\begin{array}{cc}
0 & i\\
-i & 0
\end{array}\right),\ \sigma_3=\left(\begin{array}{cc}
1 & 0\\
0 & -1
\end{array}\right).$$
These matrices are called Pauli spin matrices.

By the action of ${\rm SL}(2,\C)$ on $\H^3$, there exists a map $F$ from an open set
in $M$ to ${\rm SL}(2,\C)$ such that
$$F(\sigma_i)=F\sigma_iF^*=e^i,\ i=0,1,2,3.$$
This map represents the local moving frame associated to the immersion $f$. Let
$\Omega:=F^{-1}dF\in{\rm\mathfrak sl}(2,\C)$. The Gauss and Codazzi equations are
equivalent to
$$d\Omega+{\frac{1}{2}}[\Omega\wedge\Omega]=0,$$
which is the null curvature condition of the Maurer-Cartan (connection) form $\Omega$.

It is well known \cite[for example]{Bo} that every surface with constant mean
curvature (in $\E^3$, $\S^3$, $\H^3$) admits conformal (isothermal) coordinates,
$z=x+iy$, so that
$$I=ds^2=df\otimes df=e^{2u}dz\otimes d\bar{z}.$$

So, we may rewrite $f: D \longrightarrow \H^3$ (by abuse of notation), with $D \in \C$
open, simply connected, and containing the origin.

Thus, $<f_z,f_z>=<f_{\bar z},f_{\bar z}>=0$, $<f_z,f_{\bar z}>=\frac{1}{2} {e^{2u}}$,
where $f_z=\dis\frac{1}{2}(f_x-if_y)$, $f_{\bar z}=\dis\frac{1}{2}(f_x+if_y)$. We also
have
$$<f_z,N>=<f_{\bar z},N>=0,\ <N,N>=1.$$
The form $Qdz^2:=<f_{zz},N>dz^2$ is called Hopf differential.

It is also well known \cite{Bo} that an immersion $f$ has constant mean curvature if
and only if the Hopf differential is holomorphic. The second fundamental form is
defined as
$$II=-<df,dN>=ldx^2+2mdxdy+ndy^2.$$ Then
\begin{eqnarray*}
<f_{zz},N>&=&\frac{1}{4}(l-n-2im)=Q,\\
<f_{z{\bar z}},N>&=&\frac{1}{4}(l+n)=\frac{1}{2}He^{2u},
\end{eqnarray*}
where $N\equiv e_3$ represents the usual Gauss map (unit normal vector field on $M$).
The Maurer-Cartan form $\Omega$ can be written as
$$\Omega=Adz+Bd{\bar z},$$
where
$$A=\left(\begin{array}{cc}
\frac{1}{2}u_z & \frac{1}{2}e^u(1+H)\\
-e^{-u}Q & -\frac{1}{2}u_z
\end{array}\right),\
B=\left(\begin{array}{cc}
-\frac{1}{2}u_{\bar z} & e^{-u}{\overline Q}\\
\frac{1}{2}e^u(1-H) & \frac{1}{2}u_{\bar z}
\end{array}\right).$$
The moving frame $F$ satisfies the following Lax equations
\begin{equation}
\left\{\begin{array}{ccc}
F_z&=&FA\\
F_{\bar z}&=&FB
\end{array}
\right..
\end{equation}
The compatibility condition $F_{z\bar z}=F_{{\bar z}z}$ gives
\begin{equation}
A_{\bar z}-B_z-[A,B]=0,
\end{equation}
which can be written as
\begin{equation}
\label{eq:gc} \left\{\begin{array}{c}
u_{z{\bar z}}-\dis\frac{e^{2u}}{4}(1-H^2)-{e^{-2u}}Q{\overline Q}=0,\\
Q_{\bar z}=0.
\end{array}
\right.
\end{equation}

\section{Spectral Deformations}
Let us consider an arbitrary immersion $f$, of metric conformal factor $e^{u(z, \bar
z)}$, constant mean curvature $H$ and Hopf differential $Q dz^2$. We start with the
desire to find isometric or quasi-isometric spectral deformations of this immersion.
We are mainly interested in new surfaces characterized by the triple
$(\tilde{u},\tilde{H},\tilde{Q})$, such that they satisfy the same Gauss and Codazzi
equations as the initial $(u,H,Q)$. Here, it should be remarked that our approach is
different from the one in \cite{Bo} and \cite{Qi}. More precisely, these references
consider a spectral transformation given by the complex non-zero parameter $\lambda $
such that
$$(1+H)\to{\lambda}(1+H),\ (1-H)\to{\lambda^{-1}}(1-H)$$ (see, for example,
\cite[formulas 23-25]{Qi}); this accordingly changes the matrices $A$ and $B$ of the
Lax system.

$H$ and $u$ are both real-valued, while $Q$ is complex-valued. In our opinion, the
following two spectral deformations have more geometric meaning: \

A). A {\it positive real parameter} $s$ is introduced in the second term of the
Gauss-Codazzi equation via $(1+H)\to s(1+H)$ and $(1-H)\to s^{-1} (1-H)$, so that the
new $(\tilde{u},\tilde{H},\tilde{Q})$ satisfy the same Gauss-Codazzi equation. \

B). A {\it complex parameter of modulus one}, $\theta=e^{it}$, is introduced in the
third term of the Gauss-Codazzi equation, via $Q\to {\theta^{-2}} Q$, so that the
Gauss-Codazzi equation does not change.

\

A). {\bf $s$-Spectral Deformations. Proper Deformations}

The name of spectral parameter comes from mathematical physics, where it was
interpreted as a spectral parameter in a corresponding linear problem.

\begin{definition}
We call {\it s-spectral deformation} of the constant mean curvature immersion $f$ the
effect (on the surface) of introducing the positive parameter $s$ via $(1+H)\to
s(1+H)$ and $(1-H)\to s^{-1}(1-H)$, respectively.

\

This effect depends on the geometric interpretation we give this transformation, that
is:
\begin{equation}
k(1+ H^s):=s(1+H)
\end{equation}
and
\begin{equation}
k(1- H^s): = s^{-1}(1-H),
\end{equation}
where $k$ is a nonzero real number.
\end{definition}
As a direct consequence of equations above, we obtain :
\begin{equation}
k=\frac{s(1+H)+s^{-1}(1-H)}{2}
\end{equation}
\begin{equation}
H^s=\frac{s(1+H)-s^{-1}(1-H)}{s(1+H)+s^{-1}(1-H)}
\end{equation}
\begin{theorem}
\label{th:s-spec} For any fixed positive parameter $s$, the $s$-spectral
transformation
\begin{equation}
(1+H)\to s(1+H)(=k(1+ H^s))
\end{equation}
and
\begin{equation}
(1-H)\to s^{-1}(1-H)(=k(1- H^s)),
\end{equation}
deforms an immersion $f$ of metric $e^{2u}dzd{\bar z}$, Hopf differential $ Qdz^2$ and
mean curvature $H$, into a conformal immersion, $f^s$, of metric $e^{2u^s}dz d{\bar
z}:= |k|^2 \cdot e^{2u}dzd{\bar z}$, Hopf differential $Q^s dz^2 := k \cdot Q dz^2$
and mean curvature $H^s$, as defined by the formulas above.
\end{theorem}
\begin{proof}
Note that the Gauss-Codazzi equation satisfied by $(u,H,Q)$ is equivalent to the
following Gauss-Codazzi equation satisfied by $(u^s,H^s,Q^s)$:
$$\left\{\begin{array}{c}
u^s_{z{\bar z}}-\dis\frac{e^{2 u^s}}{4}[1-(H^s)^2]-e^{-2u^s}{Q^s}{\overline
Q}^s=0,\\
{Q^s}_{\bar z}=0,
\end{array}
\right.
$$
where $u^s:=u+\ln|k|$, $H^s=\dis\frac{s(1+H)-s^{-1}(1-H)}{s(1+H)+s^{-1}(1-H)}$, and
$Q^s=kQ$.
\end{proof}
This spectral deformation of $f$ to $f^s$ may be interpreted as a substitute for
similarity transformations, since similarity does not exist in hyperbolic 3-space
$\H^3(-1)$.
\begin{remark}
Any such $s$-spectral deformation is interesting in itself; it {\it rescales} both the
metric I and the 2-form $|Q|^2 dz \cdot d{\bar z}$ (by multiplication with the same
positive constant) and so the new surface looks similar to the first one, although the
mean curvature changes.
\end{remark}

We say that the immersion $f^s$ obtained from $f$ via the spectral deformations (8)
and (9) is {\it strongly conformal} to $f$.

\begin{definition}
The $s$-spectral deformation is called {\it proper} if $s\ne 1$ {\it and} it leaves
the metric unmodified, that is, $k=1$ or $k=-1$ (see Theorem \ref{th:s-spec}).
\end{definition}

\

Note that whenever $s$ is not equal to $1$, the deformation is proper iff
$s=\dis\frac{|1-H|}{|1+H|}$, $H\neq -1, 1$.

In this work, we will use general $s$-deformations (strongly conformal deformations),
and will specify those particular instances when deformations are proper (isometric).
\

B). {\bf $\theta$-Spectral Deformation:}
\begin{definition}
We call $\theta$-{\it spectral deformation} of the constant mean curvature immersion
$f$ the effect of introducing the $\S^1$-parameter $\theta=e^{it}$ such that the Hopf
differential changes according to $Q\to{\theta}^{-2} Q.$
\end{definition}
The $\theta-$deformation does not change the metric or the mean curvature, only the
Hopf differential. It gives the well-known family of {\it associate surfaces}.
\section{The $\lambda$-Spectral Deformation }
Let us consider a simply connected Riemann surface, immersed in $\H^3$. Let the
immersion be $f$, of constant mean curvature $H$ and Hopf differential $Qdz^2$.
\begin{remark}
It is easy to see that the two types of spectral deformations have different geometric
effects on the surface. We will combine the two deformations, and introduce a
parameter that covers both spectral deformations mentioned above.
\end{remark}
\begin{definition}
\label{def:lambda} We define $\lambda = s \cdot {\theta}$, where $s>0$, and $\theta
=e^{it}$. We call $\lambda $ {\it generalized spectral parameter}. \end{definition}

For the case of an {\it isometric deformation} (that is a $\theta$-deformation while
$s=1$, or a $theta$-deformation combined with a proper $s$-deformation), the mean
curvature $H$ and Hopf differential $Q$ remain the same - up to an eventual change in
sign.

\begin{definition}
By $\lambda (= s \cdot \theta)$-spectral deformation we mean the effect of performing
both of the following deformations on the initial immersion $f$ or mean curvature $H$:

A). an $s$-deformation ($s>0$), \

B). a $\theta=e^{it}$-deformation.

\end{definition}

Note that order does not matter: since $s$-deformations are independent from
$\theta$-deformations, they commute.

Case A). gives a genuine (and strongly conformal) surface deformation in general, as
described in Theorem 3. The surface stays the same in just two cases: the trivial case
$s=1$ (identity) and the case of proper deformation ($s=\frac{|1-H|}{|1+H|}$), both
being isometries.

Case B). describes the associate family.

\

While performing a general $\lambda$-deformation, that is a $\theta$-deformation and
an $s$-deformation, keep in mind the changes described in Theorem 2. In terms of Lax
matrices, we obtain:

\begin{equation}
A(s, \theta) =\left(\begin{array}{cc}
\frac{1}{2}u_z & \frac{s}{2}e^u(1+H)\\
-e^{-u} {\theta^{-2}Q} & -\frac{1}{2}u_z
\end{array}
\right),
\end{equation}

\begin{equation}
B(s, \theta)=\left(\begin{array}{cc}
-\frac{1}{2}u_{\bar z} & e^{-u}{\theta^{2}\overline Q}\\
\frac{s^{-1}}{2}e^u(1-H) & \frac{1}{2}u_{\bar z}
\end{array}
\right).
\end{equation}

\begin{remark}

For loop group reasons, we conjugate these matrices with the $z$-independent matrix
\begin{equation} G=i \left(\begin{array}{cc}
0 & \theta^{1/2}\\
\theta^{-1/2} & 0
\end{array}
\right),
\end{equation}

and obtain the matrices

\begin{equation}
A^{\lambda} =\left(\begin{array}{cc}
-\frac{1}{2}u_z & -\theta^{-1}\cdot e^{-u}Q\\
\theta^{-1}\cdot \frac{s}{2}e^u(1+H) & \frac{1}{2}u_z
\end{array}
\right),
\end{equation}

\begin{equation}
B^{\lambda}=\left(\begin{array}{cc}
\frac{1}{2}u_{\bar z} & \theta \cdot \frac{s^{-1}}{2}e^u(1-H)\\
\theta \cdot e^{-u}{\overline Q} & -\frac{1}{2}u_{\bar z}
\end{array}
\right).
\end{equation}

Note that these conjugated matrices will satisfy the Lax system and the compatibility
condition associated to it.

\begin{remark}

While looking for the right type of spectral transformation, eventually an isometric
one, instead of our deformation, one may have been tempted to perform the traditional
one: $Q \longrightarrow \lambda^{-1} Q$, with $\lambda \in {\Bbb C}^*$, hoping to
obtain a frame $F$ with the property $F\cdot {\overline {F({\bar {\lambda}}^{-1})}}^t
= I$, which in particular would be unitary for $\lambda$ in $S^1$. Note that this type
a deformation does not lead to such a frame.

Also, if we made such a choice, the off-diagonal terms that contain $H+1$ and $H-1$
would destroy the hope for a $\rm su(2)$-valued Maurer-Cartan form.

\ The $\lambda$ deformation we just introduced is convenient, in the sense that the
Maurer-Cartan form becomes a $\rm su(2)$-valued form for a specific real value $s_0$
of the parameter $s$, and all values of $\theta$ in $S^1$.
\end{remark}

\

Via our $\lambda$ deformation, the frame $F$ changes to $F^{\lambda}$ (which can be
considered fixed at a point $p\in M$). The Lax system
$$\left\{\begin{array}{ccc}
F^\lambda_z&=&F^\lambda A^\lambda\\
F^\lambda_{\bar z}&=&F^\lambda B^\lambda
\end{array}
\right.$$ can be also written as
\begin{equation}
\label{mceqn} (F^\lambda)^{-1}dF^\lambda=\Omega^\lambda,
\end{equation}
so that the Maurer-Cartan form $\Omega^\lambda$ writes
\begin{equation}
\Omega^\lambda=A^\lambda dz+B^\lambda d{\bar z}
\end{equation}

\

A solution $F^\lambda$ of the above equation, together with the initial condition
$F^\lambda(0,0,{\lambda})=I$, in a simply connected domain $D$, $F^\lambda:
D\longrightarrow\Lambda^s{\rm SL}(2,\C)$, is called extended frame corresponding to
the spectral deformations $f\longmapsto f^s$, and $Q\longmapsto\theta^{-2}Q$.

\

Here, $\Lambda^s{\rm SL}(2,\C)$ represents the ``twisted'' loop group over ${\rm
SL}(2,\C)$ given by the automorphism
$$\sigma: g\longmapsto({\rm Ad}\sigma_3)(g),$$
$$\sigma_3=\left(\begin{array}{cc}
1 & 0\\
0 & -1
\end{array}\right),$$

$\Lambda^s{\rm SL}(2,\C):=\{g: C_s\longrightarrow{\rm SL}(2, \C)|
g(-\lambda)=\sigma(g(\lambda))\},$ where $s$ is the absolute value of the parameter
$\lambda$ and $C_s$ is the circle of center O and radius $s$ in the complex plane.

Note these maps could be also written as $g_s(\theta): {S_1}\longrightarrow{\rm SL}(2,
\Bbb C),$ with $s$ fixed, real and positive, and the property $g_s
(-\theta)=\sigma(g_s( \theta))$, and that there is no significant difference between
these loop groups and the usual loop group considered in \cite{DoHa}.

\

We will denote $\Lambda{\rm SL}(2,\C):=\{g: {S^1}\longrightarrow{\rm SL}(2, \C)|
g(-\theta)=\sigma(g(\theta))\},$

\end{remark}

It is customary to denote by $\Lambda_-^*{\rm SL(2,\C)}$ the set of all maps of
$\Lambda{\rm SL(2,\C)}$ that can be holomorphically extended outside the disk enclosed
by the circle, and equal to identity at infinity. Also, $\Lambda_+{\rm SL(2,\C)}$
stands for those maps that can be holomorphically extended inside the same disk.

 Similar notations are used for $\Lambda{\rm SU}(2)$.

 In order to make such a loop group into a complete Banach Lie group,
we consider the same $H^p$-norm for $p>\frac{1}{2}$ as used in \cite{DoHa}. Elements
of this loop group are matrices with off-diagonal entries that are  odd in $\theta$
and diagonal entries that are even in $\theta$. We view the elements as {\it formal
series in $\theta$}.

\

Whenever we use loop group factorizations, we will always split inside the loop group
$\Lambda{\rm SL(2,\C)}$. The reason why we use loop group factorizations is related to
the methods of constructing surfaces starting from the generalized Weierstrass
representation formula. Such a method was first presented in \cite{DoHa}.

\begin{theorem} For any associated family of CMC surfaces with given frame $F=F(\theta)$, $\theta \in S^1$,
 and mean curvature $|H|>1$, there exists a certain s-deformation, for some $s=s_0$, that generates a unitary
frame $\tilde F = \tilde F(\theta) \in \Lambda{\rm SU}(2)$. The unitary frame $\tilde
F$ represents the lift of a harmonic map $\tilde N$ in $S^2$.

\end{theorem}

\begin{proof}
It is easy to see that choosing $s = s_0 := \sqrt{\frac{H-1}{H+1}}$ gives the only
deformation that makes (changes) the Maurer-Cartan $\Omega$ into an $\rm su(2)$-valued
form $\tilde \Omega$.

 Remark that as $\lambda$ we approaches $\lambda_0 = s_0 \cdot \theta$, the mean curvature will
go to infinity, and this particular deformation degenerates. From the Gauss-Codazzi
equations, it follows that there exists a map $\tilde F$ from $D$ to $SU(2)$ such that
${\tilde F}^{-1} d{\tilde F} = \tilde \Omega$. The harmonic map $\tilde N$ represents
the natural projection of the frame $\tilde F$ to $S^2$.
\end{proof}

\begin{definition}
We call $\tilde F$ the {\it adjusted frame} of $F$ and the form ${\tilde F}^{-1}
d{\tilde F}$ the {\it adjusted Maurer-Cartan form}.

\end{definition}

Hence, the explicit form of the adjusted Maurer-Cartan is

\begin{equation}
{\tilde \Omega} =\left(\begin{array}{cc}
-\frac{1}{2}u_z & -\theta^{-1}\cdot e^{-u}Q\\
\theta^{-1}\cdot \frac{1}{2}e^u {\sqrt{H^2-1}} & \frac{1}{2}u_z
\end{array}
\right) dz + \left(\begin{array}{cc}
\frac{1}{2}u_{\bar z} & - \theta \cdot \frac{1}{2}e^u {\sqrt{H^2-1}}\\
\theta \cdot e^{-u}{\overline Q} & -\frac{1}{2}u_{\bar z}
\end{array}
\right) d{\bar z}.
\end{equation}

\section{Weierstrass Type Representation Formula for CMC Surfaces in $\H^3$}

Let $M$ be any simply connected Riemann surface immersed in $\H^3$, via immersion $f$,
corresponding to the moving frame $F$.

It is well-known that for every local framing $F$ and connection form
$\Omega:=F^{-1}dF$, we have the identity (Maurer-Cartan equation):
$$d\Omega+\frac{1}{2}[\Omega\wedge\Omega]=0.$$

An arbitrary $\lambda \in {\Bbb C}^*$ deformation transforms $\Omega $ into
$\Omega^\lambda := (F^\lambda)^{-1}dF^\lambda=A^\lambda dz+B^\lambda d{\bar z},$ which
can be also written as
$$\Omega^\lambda =\Omega^{'}_1dz+\Omega_0+\Omega^{''}_1d{\bar
z},$$ where $\Omega_0=\Omega^{'}_0 dz+\Omega^{''}_0 d{\bar z}$. Here $\Omega_0$ is, as
usual, a one form with values on the diagonal elements of $\mathfrak{sl}(2,\C)$, and
the rest of the terms are off-diagonal. Hence,
\begin{eqnarray*}
\Omega^{'}_0=\left(\begin{array}{cc}
-\frac{1}{2}u_z & 0\\
0 & \frac{1}{2}u_z
\end{array}\right),\\
\Omega^{''}_0=\left(\begin{array}{cc}
\frac{1}{2}u_{\bar z} & 0\\
0 & -\frac{1}{2}u_{\bar z}
\end{array}\right),\\
\Omega^{'}_1=\left(\begin{array}{cc}
0 & -\theta^{-1}\cdot e^{-u}Q\\
\theta^{-1}\cdot \frac{s}{2}e^u(1+H) & 0
\end{array}\right),\\
\Omega^{''}_1=\left(\begin{array}{cc}
0 & \theta \cdot \frac{s^{-1}}{2}e^u(1-H)\\
\theta \cdot e^{-u}{\overline Q} & 0
\end{array}\right).
\end{eqnarray*}

Here $\lambda = s \cdot \theta$ and $\theta = e^{it}$, as usual.

Let us consider the associate form $\tilde \Omega$ which is a
$\mathfrak{su}(2)$-valued 1-form and hence decomposed via Cartan decomposition, as
$\mathfrak{su}(2)={\mathfrak k}\oplus{\mathfrak p}$, where ${\mathfrak k}$ is the
diagonal part, and ${\mathfrak p}$ is the off-diagonal one. Thus, the form $\tilde
\Omega$ writes $\tilde \Omega = \tilde \Omega_0+\tilde \Omega_1$. Let $\tilde
\Omega_1=\tilde \Omega^{'}_1 + \tilde \Omega^{''}_1$ be the splitting into $(1,0)$
and, respectively, $(0,1)$-forms. We compare matrices $\tilde \Omega^{'}_1$ and
$\tilde \Omega^{''}_1$ to their correspondents from $\Omega^{\lambda}$, namely
$\Omega^{'}_1$ and $\Omega^{''}_1$. A straight-forward computation leads us to the
following result:

\begin{theorem}
Let $\tilde N: D \longrightarrow S^2$ be a non-conformal harmonic map with lift
$\tilde F = \tilde F (\theta): D \longrightarrow {\Lambda{\rm SU}(2)}$, where $D$ is a
simply connected domain as before.

Let ${\tilde \Omega(\theta)} = {\tilde F}^{-1} d{\tilde F} = \tilde \Omega^{'}_1 +
\tilde \Omega_0 + \tilde \Omega^{''}_1$. Let $a>0$ be an arbitrary real constant, and
let

$$\beta^{'}_1 (a)= {\frac{1}{4}}\cdot(a-1)\cdot (\sigma_0 - \sigma_3) {\tilde \Omega^{'}_1 }
(\sigma_0 + \sigma_3) dz,$$

respectively

$$\beta^{''}_1 (a)= {\frac{1}{4}}\cdot(a^{-1}-1)\cdot (\sigma_0 + \sigma_3) {\tilde \Omega^{''}_1 }
(\sigma_0 - \sigma_3) d{\bar z} $$

Let $$\Omega = \Omega(a, \theta) :={\tilde \Omega (\theta)} + {\beta^{'}_1(a)} +
{\beta^{''}_1(a)}.$$

Then we have the following:

i).$$d\Omega (a, \theta)+\frac{1}{2}[\Omega (a, \theta)\wedge\Omega (a, \theta)]=0.$$

ii).If $F$ is a $\rm SL(2, \C)$-valued solution of $\Omega = F^{-1}dF$, then $f = F
\cdot F^*$ is a conformal immersion with isolated singularities and constant mean
curvature $H = \frac{a^2 + 1}{a^2 -1}$.

\end{theorem}

\begin{proof}

One may also see Theorem 4.4, \cite{Qi}, for the construction of a form that is
similar to $\Omega$.

Let us now assume that $\tilde \Omega (\theta) = \tilde F^{-1} d {\tilde F}$ is of the
form given by equation (17). We are looking for a parameter $\lambda \in \C^*$ such
that $\Omega^{\lambda}$ coincides with $\Omega (a, \theta)$. By direct computation, we
obtain that $s=|\lambda|$ must satisfy the relation $s=a\cdot
{\sqrt{\frac{H-1}{H+1}}}$. Therefore, $\lambda = s \cdot \theta \in C_s$, where $s$ is
uniquely determined by $a$, from the above formula. The frame $F^{\lambda}$ is a
solution to $\Omega = F^{-1}dF$, and its corresponding immersion is $f^\lambda =
F^\lambda \cdot {F^\lambda}^*$, with mean curvature $H^s$. We substitute $s$ in the
simplified formula (7) and obtain ${H^s}=\frac{a^2 - 1}{a^2 + 1}$. This proves ii).
\end{proof}

\begin{remark}

Assuming that the initial Maurer-Cartan form $\tilde \Omega$ corresponds to the
adjusted frame $\tilde F$ of a certain family of CMC surfaces with initial mean
curvature $H$, let us compute the matrices explicitely. We obtain

$$\beta^{'}_1 = (a-1) \left(\begin{array}{cc}
0 & \qquad 0\\
\theta^{-1} \cdot \frac{1}{2} e^u {\sqrt{H^2-1}}& \qquad 0
\end{array}\right)dz,$$

respectively

$$\beta^{''}_1= (a^{-1}-1) \left(\begin{array}{cc}
0 & \qquad -\theta \cdot \frac{1}{2} e^u {\sqrt{H^2-1}}\\
0 & \qquad 0
\end{array}\right)d{\bar z} $$

The sum of these matrices is $\Lambda {\rm sl}(2,\C)$-valued. Note that the first
defines a $(1,0)$ form in $\theta^{-1}$, while the other one is a $(0,1)$ form in
$\theta$. These expressions will be of use in the next section.

Note that while $\Omega$ is not ${\rm su}(2)$-valued, $\beta^{'}_1+\beta^{''}_1$
measures its `defect' from ${\rm su}(2)$. In a sense, this measures the `defect' of
the Gauss map $N$ from being a harmonic map in the symmetric space $S^2$.
\end{remark}

Let us now denote

\begin{equation}
G(\theta) ={F^{\lambda}} \cdot \tilde F^{-1}(\theta)
\end{equation}

A very important remark is that $G \cdot G^* = F \cdot F^* = f^{\lambda}$. We view the
matrix $G$ exclusively as a function of $\theta$. It will play a significant role in
Section 7.

\section{Normalized potentials}

The notion of normalized potential was introduced in the most general case - for
harmonic maps in symmetric spaces, and their extended frames - \cite{DoPeWu}. Next,
\cite{DoHa} and \cite{Wu} gave the expression and computation of this potential in
particular for the case of constant mean curvature surfaces in Euclidean space. We
recall the following adaptation of (see \cite{Wu}):

\begin{theorem}
Let ${\tilde N}: D \rightarrow S^2$ be a harmonic map based at identity, and $\tilde
F(\theta ):D\rightarrow \Lambda{\rm SU}(2, \C)$ an extended frame corresponding to it.
Then there exists a discrete subset $S$ of $D-{0}$ such that for any $z \in {D-S}$ we
have ${\tilde F}(z,\theta)={\tilde F}_-(z,\theta)\cdot {\tilde F}_+(z,\theta)$, with
${\tilde F}_-(\theta) \in \Lambda_-^*{\rm SL}(2,\C)$ and ${\tilde F}_+(\theta) \in
\Lambda_+{\rm SL}(2,\C)$. Remark that here minus and plus refer to the power series in
$\theta$. The form $P(z) = {\tilde F}_-^{-1} {d{\tilde F}_-} \theta$ is a meromorphic
(1,0) - form on $D$, with poles in $S$. This form is called {\it meromorphic
potential} or {\it normalized potential}.

\

Conversely, any such harmonic map $\tilde N$ can be constructed from a meromorphic
potential by integration, obtaining first ${\tilde F}_- : D-S \in {\Lambda_-^*{\rm
SL}(2,\C)}$ where the discrete subset $S$ consists of poles of $P$ and then obtaining
an extended frame ${\tilde F}$ of $f$ via the Iwasawa factorization $\Lambda{\rm
SL}(2,\C)=\Lambda{\rm SU}(2)\cdot \Lambda_+^B{\rm SL}(2,\C)$, ${\tilde F}_-={\tilde
F}\cdot {\tilde F}_+^{-1}$.

\end{theorem}

For details on the Iwasawa factorization, one may consult \cite{Pr} and \cite{DoHa}.
Note that this type of decomposition may be done in minus-plus form or in plus-minus
form (with different, unique factors).

The above stated theorem has the following important consequence:

\begin{theorem}
The normalized potential corresponding to constant mean surfaces in the hyperbolic
space is identical to the one corresponding to their Euclidean correspondents.
\end{theorem}

\begin{proof} The Lawson correspondence is performed via the same harmonic maps. More
precisely, the harmonic maps that represent Gauss maps for the Euclidean CMC surfaces
correspond to the adjusted Gauss maps of their hyperbolic counterparts. Via the above
theorem, it becomes natural that the normalized potential in the two cases is the
same.

\end{proof}

In view of the above theorem and using formula (3.24) of \cite{Wu}, we deduce the
normalized potential corresponding to a CMC surface with $|H|>1$ in hyperbolic space,
of Hopf differential $Q dz^2$ and metric factor $u(z, \bar z)$, as

\

$P = \left(\begin{array}{cc}
0 & -e^{-2h(z)+h(0)}Q\\
\frac{1}{2}e^{2 h(z)-h(0)}\sqrt{H^2-1} & 0
\end{array}\right)$

\

where $h(z):=u(z,0)$ is the {\it holomorphic part} of $u(z,\bar z)$.

Remark that $\theta^{-1} P$ can be deduced directly from the form $\tilde
\Omega^{'}_1$, as we had expected.

\

Note that we did not use Lawson's correspondence in order to obtain this result. In
some other words, we did not `cheat', by replacing some Euclidean mean curvature $c$
with its hyperbolic correspondent $\sqrt{c^2-1}$. The normalized potential $P$ that we
arrived at simply came as a byproduct of our loop group techniques!

\

In the spirit of Wu, the holomorphic part $e^{2 h(z)}$ of the conformal factor $e^{2
u(z,\bar z)}$ in the induced metric $e^{2u} dz d{\bar z}$ on a CMC immersion $f$ in
the hyperbolic space is meromorphic on the entire domain $D$. This meromorphic
function and the Hopf differential uniquely determine the induced metric and the
surface, up to spectral deformations.

The Weierstrass type data (potential) is the ``genetic material" (like the classical
Weierstrass representation formula for minimal surfaces in $\E^3$) for surface
construction.

\section{Constructing CMC surfaces in $\H^3$}

This represents a DPW type of algorithm to be used in constructing CMC surfaces in
$\H^3$:

i). Start from a normalized potential $P$.

Solve the initial value problem $$\theta^{-1} P=\tilde F_-^{-1} d \tilde F_-$$ with
$\tilde F_-(z=0,\theta) = I$. \

ii). Given the solution $\tilde F_-$ from i), perform the Iwasawa decomposition in
$\Lambda{\rm SL}(2,\C)=\Lambda{\rm SU}(2)\cdot \Lambda_+^B{\rm SL}(2,\C)$, $$\tilde
F_-=\tilde F\cdot \tilde F_+^{-1}$$ in order to obtain the $\Lambda {\rm SU(2)}$ -
valued extended frame $\tilde F$.

Denote by $\tilde F_0$ be the coefficient of $\theta^0$ in the $\theta $-series
expansion of $\tilde F$. \

iii). Write $\tilde\Omega =\tilde F^{-1} d\tilde F$ and name it adjusted Maurer-Cartan
form. Separate its components $\tilde \Omega^{'}_1$ and $\tilde \Omega^{''}_1$. Pick a
value $a>0$. Apply formulas of Theorem 4 in order to obtain $\beta^{'}_1$ and
$\beta^{''}_1$. \

iv). Using $\beta^{'}_1$ from iii), solve the initial value problem $$ G_-^{-1} dG_- =
{\tilde F}_0 \cdot \beta^{'}_1 \cdot {\tilde F}_0^{-1},$$ with $G_-(z=0,\theta)=I$. \

v). Using $\beta^{''}_1$ from iii), solve the initial value problem $$ G_+^{-1} dG_+=
{\tilde F}_0 \cdot \beta^{''}_1 \cdot {\tilde F}_0^{-1},$$ with $G_+(z=0,\theta)=I$. \

vi). Compute $$L:=G_-^{-1}\cdot G_+,$$ where $G_-$ is the solution of iv), and $G_+$
is the solution of v).

Then split again, $$L = p_+ \cdot p_-^{-1}.$$\

Let $$G: = G_-\cdot p_+.$$ \

vii). $f := G G^*$ will be an associate family of immersions of mean curvature
$H=\frac{a^2+1}{a^2-1}$, corresponding to the normalized potential (Weierstrass
representation) $P$. (Note that there is such a family for each value of $a>0$).

\

The extended frame of this immersion is the $\Lambda {\rm SL}(2, \C)$-valued frame $F
:= G\cdot \tilde F$, and $f=F F^*$. The associate frame is $\tilde F$, which is
$SU(2)$-valued. The Gauss map $N$ of this immersion is not harmonic, but the adjusted
Gauss map $\tilde N$ is harmonic, and it represents the natural projection on $S^2$ of
the associate frame $\tilde F$.

\begin{proof}
Steps i) and ii) represent a parallel to the standard DPW procedure for CMC surfaces
in Euclidean space. They come as a direct consequence of Theorems 5 and 6.

Step iii) is justified by Theorem 4.

For the rest of the steps, let us note the consequences of defining $F := G\cdot
\tilde F$. A direct consequence is that

$$G^{-1} dG = {\tilde F} [\beta^{'}_1 + \beta^{''}_1] {\tilde F}^{-1}.$$

Iwasawa factorization for $G$ (both ways) give $G_-$ and $G_+$, uniquely, such that

$G= G_-\cdot p_+$ and $G=G_+ \cdot p_-$. After replacing these expressions of $G$ into
the formula for $G^{-1} dG$ and differentiating, we compare the terms corresponding to
negative respectively positive powers of $\theta$. We obtain

$$ G_-^{-1} dG_- = {\tilde F}_0 \cdot \beta^{'}_1 \cdot {\tilde F}_0^{-1},$$

$$ G_+^{-1} dG_+ = {\tilde F}_0 \cdot \beta^{''}_1 \cdot {\tilde F}_0^{-1},$$

where $\tilde F_0$ represents the constant matrix in the $\theta$-expansion of $\tilde
F$.

The rest of the steps are clear, and represent a standard technique of regaining a
matrix from the first factor of its Iwasawa factorization. Since $F := G\cdot \tilde
F$, and $\tilde F$ is unitary, the CMC immersion is obtained as $f=F F^*=G G^*$. The
mean curvature of the associated family of $f$ is $\frac{a^2+1}{a^2-1}$, as a
consequence of Theorem 4.

\end{proof}

\section{Open Problems}
Besides its direct computational applications, the above representation formula will
hopefully lead to a better understanding of global period problems for CMC immersions
in hyperbolic space, as well as symmetries and singularities.

Acknowlegement: The author is grateful to Prof. J. Dorfmeister for his most useful
comments and suggestions regarding this work.

\

\end{document}